\documentclass[onecolumn]{autart}
\usepackage{eufrak,mathrsfs,amsmath,longtable}
\input{amssym}

\newcommand{\p}{\rm \partial}

\newcommand{\EE}{\mathscr E}

\def\RR{\mathop{\Bbb R}}

\renewcommand{\theequation}{\arabic{section}.\arabic{equation}}
\begin{document}
\begin{frontmatter}
\title{Symmetry group classification for general\\
Burgers' equation}
\thanks[footnoteinfo]{ Corresponding author: Tel. +9821-73913426.
Fax +9821-77240472.}
\author[]{M. Nadjafikhah\thanksref{footnoteinfo}}\ead{m\_nadjafikhah@iust.ac.ir},
\author[]{R. Bakhshandeh-Chamazkoti}\ead{r\_bakhshandeh@iust.ac.ir},
\address{School of Mathematics, Iran University of Science and Technology, Narmak, Tehran 1684613114, Iran.}
\renewcommand{\sectionmark}[1]{}
\begin{abstract}
The present paper solves  the problem of the group classification
of the general Burgers' equation $u_t=f(x,u)u_x^2+g(x,u)u_{xx}$,
where $f$  and $g$ are arbitrary smooth functions of the variable
$x$ and $u$, by using  Lie method. The paper is one of the few
applications of an algebraic approach to the problem of group
classification: the method of preliminary group classification. A
number of new interesting nonlinear invariant models which have
nontrivial invariance algebras are obtained. The result of the
work is a wide class of equations summarized in table form.
%
\end{abstract}
\end{frontmatter}
\section{Introduction}
It is well known that the symmetry group method plays an important
role in the analysis of differential equations.   The history of
group classification methods goes back to Sophus Lie. The first
paper on this subject is \cite{[1]}, where Lie proves that a
linear two-dimensional second-order PDE may admit at most a
three-parameter invariance group (apart from the trivial infinite
parameter symmetry group, which is due to linearity).  He computed
the maximal invariance group of the one-dimensional heat
conductivity equation and utilized this symmetry to construct its
explicit solutions. Saying it the modern way, he performed
symmetry reduction of the heat equation. Nowadays symmetry
reduction is one of the most powerful tools for solving nonlinear
partial differential equations (PDEs). Recently, there have been
several generalizations of the classical Lie group method for
symmetry reductions. Ovsiannikov \cite{[2]} developed the method
of partially invariant solutions. His approach is based on the
concept of an equivalence group, which is a Lie transformation
group acting in the extended space of independent variables,
functions and their derivatives, and preserving the class of
partial differential equations under study
\cite{[10],[11],[15]}.\newline
$~~~~~$The investigation of the exact solutions plays an important
role in the study of nonlinear physical systems. A wealth of
methods have been developed to find these exact physically
significant solutions of a PDE though it is rather difficult. Some
of the most important methods are the inverse scattering method
\cite{[3]}, Darboux and B$\rm\ddot{a}$cklund transformations \cite{[4]}, Hirota
bilinear method \cite{[4],[5]}, Lie symmetry analysis \cite{[6],[7],[8]}, etc.
The paper \cite{[9]}, based on the Lie group method,
is investigated a very famous and important equation, which is the
general Burgers' equation as the form
\begin{eqnarray}
u_t=au_x^2+bu_{xx},\label{eq:1}
\end{eqnarray}
where $u=u(x,t)$ is the unknown real function, $a,b\in\RR$ and
$ab\neq0$. In the present paper, we consider the general Burgers'
equation as the form
\begin{eqnarray}
u_t=f(x,u)u_x^2+g(x,u)u_{xx},\label{eq:2}
\end{eqnarray}
where $u=u(x,t)$ is the unknown real function, $f$  and $g$ are
arbitrary smooth functions of the variable $x$ and $u$. Eq.
(\ref{eq:2}) represents the Burgers' equation combining both
dissipative and nonlinear effects, therefore appears in a wide
variety of physical applications. So it is important to lucubrate
the exact explicit solutions and similarity reductions for this
equation \cite{[13],[14]}. Here, we get the preliminary group
classification of Eq. (\ref{eq:2}) by means of Lie point symmetry,
and the constructed optimal systems of subalgebras. The knowledge
of the optimal system of subalgebras gives the possibility of
constructing the optimal system of solutions \cite{[2],[12],[16]} and
permits the generation of new solutions starting from invariant or
non-invariant solutions.
%
%
\section{Symmetry Methods}
Let a partial differential equation contains  $p$ dependent
variables and $q$ independent variables. The one-parameter Lie
group of transformations
\begin{eqnarray}
\overline{x}_i=x_i+\epsilon\xi^i(x,u)+O(\epsilon^2);\hspace{1cm}
\overline{u}_{\alpha}=u_{\alpha}+
\epsilon\varphi^{\alpha}(x,u)+O(\epsilon^2),\label{eq:3}
\end{eqnarray}
where
$\displaystyle{\xi^i=\frac{{\p}\overline{x}_i}{{\p}\epsilon}\big|_{\epsilon=0}}$,
 $i=1,\ldots,p$, and
 $\displaystyle{\varphi^{\alpha}=\frac{{\p}\overline{u}_{\alpha}}{{\p}\epsilon}\big|_{\epsilon=0}}$,
 $\alpha=1,\ldots,q$, are given.
The action of the Lie group can be recovered from that of its
associated infinitesimal generators. we consider general vector
field
\begin{eqnarray}
X=\sum_{i=1}^p\xi^i(x,u)\frac{\p}{{\p}x_i}+
\sum_{\alpha=1}^q\varphi^{\alpha}(x,u)\frac{\p}{{\p}u^{\alpha}}.\label{eq:4}
\end{eqnarray}
on the space of independent and dependent variables. Therefore,
the characteristic of the vector field $X$ given by (\ref{eq:4})
is the function
\begin{eqnarray}
Q^{\alpha}(x,u^{(1)})=\varphi^{\alpha}(x,u)-
\sum_{i=1}^p\xi^i(x,u)\frac{{\p}u^{\alpha}}{{\p}x_i},\;\;\;
\alpha=1,\ldots,q.\label{eq:5}
\end{eqnarray}
The second prolongation of the infinitesimal operator
\begin{eqnarray}
X=\xi^1(x,t,u)\frac{\p}{{\p}x}+
\xi^2(x,t,u)\frac{\p}{{\p}t}+\varphi(x,t,u)\frac{\p}{{\p}u}~.\label{eq:8}
\end{eqnarray}
obtained via the following prolongation formulas:
\begin{eqnarray}\nonumber
X^{(2)}=X+\varphi^x\frac{\p}{{\p}u_x}+\varphi^t\frac{\p}{{\p}u_t}+\varphi^{xt}\frac{\p}{{\p}u_{xt}}+
\varphi^{xx}\frac{\p}{{\p}u_{xx}}\;,\hspace{-4.5cm}\label{eq:9}
\end{eqnarray}
the coefficients are obtained by
\begin{eqnarray}
\varphi^{\iota}=D_{\iota}Q+\xi^1u_{x\iota}+\xi^2
u_{t\iota},\hspace{1cm}
\varphi^{\iota\jmath}=D_{\iota}(D_{\jmath}Q) +\xi^1
u_{x\iota\jmath}+\xi^2u_{t\iota\jmath},\label{eq:9}
\end{eqnarray}
where $Q=\varphi-\xi^1u_x-\xi^2u_t$ is the characteristic of the
vector field $X$ given by (\ref{eq:5}). For instance
\begin{eqnarray}
\varphi^x&=&D_x\varphi-u_xD_x\xi^1-u_tD_x\xi^2,\\\label{eq:23}
\varphi^t&=&D_t\varphi-u_xD_t\xi^1-u_tD_t\xi^2,\\\label{eq:24}
\varphi^{xx}&=&D_x\varphi^x-u_{xx}D_x\xi^1-u_{xt}D_x\xi^2,\\\label{eq:25}
\varphi^{xt}&=&D_t\varphi^x-u_{xx}D_t\xi^1-u_{xt}D_t\xi^2,\label{eq:26}
\end{eqnarray}
where the operators $D_x$ and $D_t$ denote the total derivatives
with respect to $x$ and $t$:
\begin{eqnarray}\nonumber
D_x=\frac{\p}{{\p}x}+u_x\frac{\p}{{\p}u}+u_{xx}\frac{\p}{{\p}u_x}+
u_{xt}\frac{\p}{{\p}u_t}+\ldots,\hspace{1cm}
D_t=\frac{\p}{{\p}t}+u_t\frac{\p}{{\p}u}+u_{tt}\frac{\p}{{\p}u_t}+
u_{tx}\frac{\p}{{\p}u_x}+\ldots,
\end{eqnarray}
By the theorem 6.5. in
\cite{[23]},
$\displaystyle{X^{(2)}[u_t-f(x,u)u_x^2-g(x,u)u_{xx}]\big|_{(1.2)}=0}$.
Since
$$X^{(2)}[u_t-f(x,u)u_x^2-g(x,u)u_{xx}]=\varphi^t-(f_x\xi^1+f_u\varphi)u_x^2
-(g_x\xi^1+g_u\varphi)u_{xx}-2f\varphi^xu_x-g\varphi^{xx},$$
therefore we obtain the following determining function:
\begin{eqnarray}
[\varphi^t-(f_x\xi^1+f_u\varphi)u_x^2
-(g_x\xi^1+g_u\varphi)u_{xx}-2f\varphi^xu_x-g\varphi^{xx}]\big|_{(1.2)}=0.\label{eq:11}
\end{eqnarray}
In the case of arbitrary $f(x,u)$ and $g(x,u)$ it follows that
\begin{eqnarray}
\xi^1=\varphi=\varphi^x=\varphi^t=\varphi^{xx}=0,
\end{eqnarray}
or
\begin{eqnarray}
\xi^1=\varphi=0,\;\;\;\;\xi^2=C.
\end{eqnarray}
Therefore, for arbitrary $f(x,u)$ and $g(x,u)$ Eq. (\ref{eq:1})
admits the one-dimensional Lie algebra ${\goth g}_1$, with the
basis
\begin{eqnarray}
X_2=\frac{\p}{{\p}t}.\label{eq:21}
\end{eqnarray}
${\goth g}_1$ is called the principle Lie algebra for Eq.
(\ref{eq:1}). So, it is remained to specify the coefficient $f$
and $g$ such that Eq. (\ref{eq:1}) admits an extension of the
principal algebra ${\goth g}_1$. Usually, the group classification
is obtained by inspecting the determining equation. But in our
case the complete solution of the determining equation
(\ref{eq:11}) is a wasteful venture. Therefore, we don't solve the
determining equation but, instead we obtain a partial group
classification of Eq. (\ref{eq:1}) via the so-called method of
preliminary group classification. This method was applied when an
equivalence group is generated by a finite-dimensional Lie algebra
${\goth g}_{\EE}$. The essential part of the method is the
classification of all nonsimilar subalgebras of ${\goth g}_{\EE}$.
Actually, the application of the method is simple and effective
when the classification is based on finite-dimensional equivalence
algebra ${\goth g}_{\EE}$.
\section{Equivalence transformations}
An equivalence transformation is a nondegenerate change of the
variables $t,x,u$ taking any equation of the form (\ref{eq:1})
into an equation of the same form, generally speaking, with
different $f(x,u)$ and $g(x,u)$. The set of all equivalence
transformations forms an equivalence group ${\EE}$. We shall find
a continuous subgroup ${\EE}_c$ of it making use of the
infinitesimal method.

We consider an operator of the group ${\EE}_c$ in the form
\begin{eqnarray}
Y=\xi^1(x,t,u)\frac{\p}{{\p}x}+
\xi^2(x,t,u)\frac{\p}{{\p}t}+\varphi(x,t,u)\frac{\p}{{\p}u}
+\mu(x,t,u,f,g)\frac{\p}{{\p}f}+\nu(x,t,u,f,g)\frac{\p}{{\p}g},\label{eq:12}
\end{eqnarray}
from the invariance conditions of Eq. (\ref{eq:1}) written as the
system:
\begin{eqnarray}
u_t&-&f(x,u)u_x^2-g(x,u)u_{xx}=0,\label{eq:13}\\
f_t&=&g_t=0,\label{eq:13'}
\end{eqnarray}
where $u$ and $f,g$ are considered as differential variables: $u$
on the space $(x,t)$ and $f,g$ on the extended space $(x,t,u)$.

The invariance conditions of the system (\ref{eq:13}) are
\begin{eqnarray}\nonumber
Y^{(2)}\big(u_t&-&f(x,u)u_x^2-g(x,u)u_{xx}\big)=0,\\
Y^{(2)}(f_t)&=& Y^{(2)}(g_t)=0,\label{eq:13"}
\end{eqnarray}
where $Y^{(2)}$ is the prolongation of the operator (\ref{eq:12}):
\begin{eqnarray}
Y^{(2)}=Y+\varphi^x\frac{\p}{{\p}u_x}+\varphi^t\frac{\p}{{\p}u_t}+\varphi^{xt}\frac{\p}{{\p}u_{xt}}
\varphi^{xx}\frac{\p}{{\p}u_{xx}}+\mu^t\frac{\p}{{\p}f_{t}}+
\nu^{t}\frac{\p}{{\p}g_t}.\label{eq:14}
\end{eqnarray}
The coefficients $\varphi^x, \varphi^t, \varphi^{xt},
\varphi^{xx}, \varphi^{tt}$ are given in (\ref{eq:9}) and the
other coefficients of (\ref{eq:14}) are obtained by applying the
prolongation procedure to differential variables $f$ and $g$ with
independent variables $(x,u)$. In view of (\ref{eq:13'}), we have
\begin{eqnarray}
\mu^t=\widetilde D_t(\mu)-f_x\widetilde D_t(\xi^1)-f_u\widetilde
D_t(\varphi), \hspace{1cm} \nu^t=\widetilde D_t(\nu)-g_x\widetilde
D_t(\xi^1)-g_u\widetilde D_t(\varphi),\label{eq:16}
\end{eqnarray}
where
$\displaystyle{\widetilde D_t=\frac{\p}{{\p}t}}$.
So, we have the following prolongation formulas:
\begin{eqnarray}
\mu^t=\mu_t-f_x\xi_t^1-f_u\varphi_t,\hspace{1cm}\nu^{t}=\nu_{t}-g_x\xi_t^1-g_u\varphi_t.\label{eq:19}
\end{eqnarray}
The invariance conditions (\ref{eq:13"}) give rise to
\begin{eqnarray}
\mu^t=\nu^t=0,\label{eq:20}
\end{eqnarray}
that is hold for every $f$ and $g$. Substituting (\ref{eq:20})
into (\ref{eq:19}), we obtain
\begin{eqnarray}
\begin{array}{ll}
\mu_t=\nu_t=0,\hspace{1cm} \xi^1_t=\varphi_t=0. \label{eq:21}
\end{array}
\end{eqnarray}
Moreover with substituting (\ref{eq:14}) into (\ref{eq:13"}) we
obtain
\begin{eqnarray}
\varphi^t-2f(x,u)u_x\varphi^x-g(x,u)\varphi^{xx}-\mu u_x^2-\nu
u_{xx}-\nu=0.\label{eq:22}
\end{eqnarray}
Substituting (\ref{eq:23}), (\ref{eq:24}) and (\ref{eq:25}) into
invariance condition (\ref{eq:22}), and introducing the relation
$u_t=f(x,u)u_x^2+g(x,u)u_{xx}$ to eliminate $u_t$ we are left with
a polynomial equation involving the various derivatives of
$u(x,t)$ whose coefficients are certain derivatives of
$\xi^1,\xi^2,\varphi,\mu$, and $\nu$. We can equate the individual
coefficients to zero, leading to the complete set of determining
equations:
\begin{eqnarray}
\xi^1&=&\xi^1(x) \\\label{eq:23} \xi^2&=&\xi^2(t)=0\\\label{eq:24}
\varphi_u&=&\xi_t^2\\\label{eq:25}
\nu&=&-g\xi_t^2+2\xi_x^1\\\label{eq:26}
\mu&=&-f\xi_t^2-f(\varphi_u-2\xi_x^1)-g\varphi_{uu}\label{eq:27}
\end{eqnarray}
so, we find that
\begin{eqnarray}\nonumber
&&\xi^1(x)=a(x),\hspace{1.5cm}\xi^2(t)=c_1t+c_2,\hspace{1cm}\varphi(x,u)=c_1u+b(x),\\
&&\hspace{6mm}\mu=-2f(c_1-a(x)),\hspace{1cm}\nu=-g(c_1-a'(x)),\label{eq:28}
\end{eqnarray}
with constants $c_1, c_2$ and two arbitrary functions $a(x)$ and
$b(x)$ such that $b''(x)=c_1-a'(x)$.\newline
$\;\;\;\;$We summarize: The class of Eq. (\ref{eq:2}) has an
infinite continuous group of equivalence transformations generated
by the following infinitesimal operators:
\begin{eqnarray}
Y=a(x)\frac{\p}{{\p}x}+
(c_1t+c_2)\frac{\p}{{\p}t}+(c_1u+b(x))\frac{\p}{{\p}u}
-2f(c_1-a(x))\frac{\p}{{\p}f}-g(c_1-a'(x))\frac{\p}{{\p}g},\label{eq:30}
\end{eqnarray}
Therefore the symmetry algebra of the Burgers' equation
(\ref{eq:2}) is spanned by the vector fields
\begin{eqnarray}
&& Y_1=t\frac{\p}{{\p}t}+u\frac{\p}{{\p}u}-2f\frac{\p}{{\p}f}-g\frac{\p}{{\p}g},
\hspace{1.9cm} Y_2=\frac{\p}{{\p}t}, \nonumber\\[-2mm]
&& \label{eq:31}\\[-2mm]
&&Y_3=a(x)\frac{\p}{{\p}x}+2fa(x)\frac{\p}{{\p}f}+ga'(x)\frac{\p}{{\p}g},\hspace{1cm}Y_4=b(x)\frac{\p}{{\p}u}.
\nonumber
\end{eqnarray}
Moreover, in the group of equivalence transformations there are
included also discrete transformations, i.e., reflections
\begin{eqnarray}
t\longrightarrow-t,\hspace{1.5cm}x\longrightarrow-x,\hspace{1.5cm}u\longrightarrow-u,\hspace{1.5cm}f\longrightarrow-f,
\hspace{1.5cm}g\longrightarrow-g.\label{eq:32}
\end{eqnarray}
\begin{table}
\caption{Commutation relations satisfied by infinitesimal
generators in (4.35) }\label{table:1}
\vspace{-0.3cm}\begin{eqnarray*}\hspace{-0.75cm}\begin{array}{llllll}
\hline
  [\,,\,]&\hspace{2.1cm}Y_1 &\hspace{2.75cm}Y_2  &\hspace{2.75cm}Y_3&\hspace{2.75cm}Y_4&\hspace{2.75cm}Y_5  \\ \hline
  Y_1  &\hspace{2.1cm} 0        &\hspace{2.75cm} 0   &\hspace{2.75cm}0  &\hspace{2.75cm}0  &\hspace{2.75cm}0     \\
  Y_2  &\hspace{2.1cm}0         &\hspace{2.75cm} 0   &\hspace{2.75cm} 0 &\hspace{2.75cm}Y_2&\hspace{2.75cm}0     \\
  Y_3  &\hspace{2.1cm}0         &\hspace{2.75cm} 0   &\hspace{2.75cm} 0 &\hspace{2.75cm}Y_3&\hspace{2.75cm}0     \\
  Y_4  &\hspace{2.1cm} 0        &\hspace{2.75cm} -Y_2&\hspace{2.75cm}-Y_3&\hspace{2.75cm}0  &\hspace{2.75cm}0     \\
  Y_5  &\hspace{2.1cm} 0        &\hspace{2.75cm} 0   &\hspace{2.75cm} 0 &\hspace{2.75cm}0  &\hspace{2.75cm}0     \\
  \hline
\end{array}\end{eqnarray*}
\end{table}
\begin{table}
\caption{Adjoint relations satisfied by infinitesimal generators
in (4.35) }\label{table:1}
\vspace{-0.3cm}\begin{eqnarray*}\hspace{-0.75cm}\begin{array}{llllll}
\hline
[\,,\,]&\hspace{2.1cm}Y_1&\hspace{2.75cm}Y_2&\hspace{2.75cm}Y_3&\hspace{2.75cm}Y_4&\hspace{2.75cm}Y_5  \\ \hline
  Y_1  &\hspace{2.1cm}Y_1&\hspace{2.75cm}Y_2&\hspace{2.75cm}Y_3&\hspace{2.75cm}Y_4&\hspace{2.75cm}Y_5     \\
  Y_2  &\hspace{2.1cm}Y_1&\hspace{2.75cm}Y_2&\hspace{2.75cm}Y_3&\hspace{2.75cm}Y_4-s\;Y_2&\hspace{2.75cm}Y_5     \\
  Y_3  &\hspace{2.1cm}Y_1&\hspace{2.75cm}Y_2&\hspace{2.75cm}Y_3&\hspace{2.75cm}Y_4-s\;Y_3 &\hspace{2.75cm}Y_5    \\
  Y_4  &\hspace{2.1cm}Y_1&\hspace{2.75cm}e^s\;Y_2&\hspace{2.75cm}e^s\;Y_3&\hspace{2.75cm}Y_4&\hspace{2.75cm}Y_5     \\
  Y_5  &\hspace{2.1cm}Y_1&\hspace{2.75cm}Y_2&\hspace{2.75cm}Y_3&\hspace{2.75cm}Y_4&\hspace{2.75cm}Y_5     \\
  \hline
\end{array}\end{eqnarray*}
\end{table}
\section{Preliminary group classification}
One can observe in many applications of group analysis that most
of extensions of the principal Lie algebra admitted by the
equation under consideration are taken from the equivalence
algebra ${\goth g}_{\EE}$. We call these extensions
$\EE$-extensions of the principal Lie algebra. The classification
of all nonequivalent equations (with respect to a given
equivalence group $G_{\EE}$,) admitting $\EE$-extensions of the
principal Lie algebra is called a preliminary group
classification. Here, $G_{\EE}$ is not necessarily the largest
equivalence group but, it can be any subgroup of the group of all
equivalence transformations.\newline
$~~~~~$So, we can take any finite-dimensional subalgebra
(desirable as large as possible) of an infinite-dimensional
algebra with basis (\ref{eq:31}) and use it for a preliminary
group classification. We select the subalgebra ${\goth g}_5$
spanned on the following operators:
\begin{eqnarray}
Y_1&=&\frac{\p}{{\p}x},\hspace{2cm}
Y_2=\frac{\p}{{\p}t},\hspace{2cm}
Y_3=\frac{\p}{{\p}u},\nonumber\\
Y_4&=&t\frac{\p}{{\p}t}+u\frac{\p}{{\p}u}-2f\frac{\p}{{\p}f}-g\frac{\p}{{\p}g},\hspace{1.9cm}
Y_5=\frac{\p}{{\p}x}+2f\frac{\p}{{\p}f}+g\frac{\p}{{\p}g}.\label{eq:33}
\end{eqnarray}
The communication relations between these vector fields is given
in Table 1. To each $s$-parameter subgroup there corresponds a
family of group invariant solutions. So, in general, it is quite
impossible to determine all possible group-invariant solutions of
a PDE. In order to minimize this search, it is useful to construct
the optimal system of solutions. It is well known that the problem
of the construction of the optimal system of solutions is
equivalent to that of the construction of the optimal system of
subalgebras \cite{[2],[12]}. Here, we will deal with the
construction of the optimal system of subalgebras of ${\goth
g}_5$.\newline
$~~~~~$Let $G$ be a Lie group, with ${\goth g}$ its Lie algebra.
Each element $T\in G$ yields inner automorphism
$T_a\longrightarrow TT_aT^{-1}$ of the group $G$. Every
automorphism of the group $G$ induces an automorphism of ${\goth
g}$. The set of all these automorphism is a Lie group called {\it
the adjoint group $G^A$}. The Lie algebra of $G^A$ is the adjoint
algebra of ${\goth g}$, defined as follows. Let us have two
infinitesimal generators $X,Y\in L$. The linear mapping ${\rm
Ad}X(Y):Y\longrightarrow[X,Y]$ is an automorphism of ${\goth g}$,
called {\it the inner derivation of ${\goth g}$}. The set of all
inner derivations ${\rm ad}X(Y)(X,Y\in{\goth g})$ together with
the Lie bracket $[{\rm Ad}X,{\rm Ad}Y]={\rm Ad}[X,Y]$ is a Lie
algebra ${\goth g}^A$ called the {\it adjoint algebra of ${\goth
g}$}. Clearly ${\goth g}^A$ is the Lie algebra of $G^A$. Two
subalgebras in ${\goth g}$ are {\it conjugate} (or {\it similar})
if there is a transformation of $G^A$ which takes one subalgebra
into the other. The collection of pairwise non-conjugate
$s$-dimensional subalgebras is the optimal system of subalgebras
of order $s$. The construction of the one-dimensional optimal
system of subalgebras can be carried out by using a global matrix
of the adjoint transformations as suggested by Ovsiannikov
\cite{[2]}. The latter problem, tends to determine a list (that is
called an {\it optimal system}) of conjugacy inequivalent
subalgebras with the property that any other subalgebra is
equivalent to a unique member of the list under some element of
the adjoint representation i.e. $\overline{{\goth h}}\,{\rm
Ad(g)}\,{\goth h}$ for some ${\rm g}$ of a considered Lie group.
Thus we will deal with the construction of the optimal system of
subalgebras of ${\goth g}_5$.\newline
$~~~~~$The adjoint action is given by the Lie series
\begin{eqnarray}
{\rm Ad}(\exp(s\,Y_i))Y_j
=Y_j-s\,[Y_i,Y_j]+\frac{s^2}{2}\,[Y_i,[Y_i,Y_j]]-\cdots,
\end{eqnarray}
where $s$ is a parameter and $i,j=1,\cdots,5$. The adjoint
representations of ${\goth g}_5$ is listed in Table 2, it consists
the separate adjoint actions of each element of ${\goth g}_5$ on
all other elements.

$~~~${\bf Theorem 4.1.} {\it An optimal system of one-dimensional
Lie subalgebras of general Burgers' equation (\ref{eq:2}) is
provided by those generated by}
\begin{eqnarray}\hspace{-0.7cm}\begin{array}{rlrl}
1)&Y^1=Y_1={\p}_t, \hspace{1mm}  &11)&Y^{11}=-Y_1+Y_5=-{\p}_t+{\p}_x+2f{\p}_f+g{\p}_g,  \\
2)&Y^2=Y_2={\p}_x, \hspace{1mm}  &12)&Y^{12}=Y_4+Y_5=t{\p}_t+{\p}_x+u{\p}_u, \\
3)&Y^3=Y_3={\p}_u, \hspace{1mm} &13)&Y^{13}=-Y_4+Y_5=-t{\p}_t+{\p}_x-u{\p}_u+4f{\p}_f+2g{\p}_g,\\
4)&Y^4=Y_4=t{\p}_t+u{\p}_u-2f{\p}_f-g{\p}_g,\hspace{1mm} &14)&Y^{14}=Y_1+Y_4+Y_5=(t+1){\p}_t+{\p}_x+u{\p}_u,\\
5)&Y^5=Y_5={\p}_x+2f{\p}_f+g{\p}_g,\hspace{1mm} &15)&
Y^{15}=-Y_1+Y_4+Y_5=(t-1){\p}_t+{\p}_x+u{\p}_u,\\
6)&Y^6=Y_1+Y_2={\p}_t+{\p}_x,\hspace{1cm} &16)&
Y^{16}=Y_1-Y_4+Y_5=(1-t){\p}_t-u{\p}_u+2f{\p}_f+g{\p}_g,\\
7)&Y^7=-Y_1+Y_2=-{\p}_t+{\p}_x,\hspace{1cm} &17)&
Y^{17}=-Y_1-Y_4+Y_5=-(1+t){\p}_t+{\p}_x-u{\p}_u+4f{\p}_f+2g{\p}_g,\\
8)&Y^8=Y_1+Y_4=(t+1){\p}_t+u{\p}_u-2f{\p}_f-g{\p}_g,&&\\
9)&Y^9=-Y_1+Y_4=(t-1){\p}_t+u{\p}_u-2f{\p}_f-g{\p}_g,\\
 10)&
Y^{10}=Y_1+Y_5={\p}_t+{\p}_x+2f{\p}_f+g{\p}_g,\\
\label{eq:16}
\end{array}\end{eqnarray}
{\bf Proof.} Let ${\goth g}_4$ is the symmetry algebra of
Eq.~(\ref{eq:2}) with adjoint representation determined in Table 2
and
\begin{eqnarray}
Y=a_1Y_1+a_2Y_2+a_3Y_3+a_4Y_4+a_5Y_5,
\end{eqnarray}
is a nonzero vector field of ${\goth g}$. We will simplify as many
of the coefficients $a_i;i=1,\ldots,5$, as possible through proper
adjoint applications on $Y$. We follow our aim in the below easy
cases:\newline
{\it Case 1:} \newline
$~~~~$ At first, assume that $a_5\neq 0$. Scaling $Y$ if
necessary, we can assume that $a_5=1$ and so we get
\begin{eqnarray}
Y=a_1Y_1+a_2Y_2+a_3Y_3+a_4Y_4+Y_5.
\end{eqnarray}
Using the table of adjoint (Table 2) , if we act on $Y$ with ${\rm
Ad}(\exp(a_2Y_2))$, the coefficient of $Y_2$ can be vanished:
\begin{eqnarray}
Y'=a_1Y_1+a_3Y_3+a_4Y_4+Y_5.
\end{eqnarray}
Then we apply ${\rm Ad}(\exp(a_3Y_3))$ on $Y'$ to cancel the
coefficient of $Y_3$:
\begin{eqnarray}
Y''=a_1Y_1+a_4Y_4+Y_5.
\end{eqnarray}
{\it Case 1a:} \\
$~~~~$ If $a_1,a_4\neq 0$ then we can make the coefficient of
$Y_1$ and $Y_4$ either $+1$ or $-1$. Thus any one-dimensional
subalgebra generated by $Y$ with $a_3, a_4\neq 0$ is equivalent to
one generated by $\pm Y_1\pm Y_4+Y_5$ which introduce parts {\it
14)}, {\it 15)}, {\it 16)} and {\it 17)} of the theorem.\newline
{\it Case 1b:} \newline
$~~~~$ For $a_1=0, a_4\neq0$ we can see that each one-dimensional
subalgebra generated by $Y$ is equivalent to one generated by $\pm
Y_4+Y_5$ which introduce parts {\it 12)} and {\it 13)} of the
theorem.\newline
{\it Case 1c:} \newline
$~~~~$ For $a_1\neq0, a_4=0$ we can see that each one-dimensional
subalgebra generated by $Y$ is equivalent to one generated by $\pm
Y_1+Y_5$ which introduce parts {\it 10)} and {\it 11)} of the
theorem.\newline
{\it Case 2:} \newline
$~~~~$ The remaining one-dimensional subalgebras are spanned by
vector fields of the form $Y$ with $a_5=0$. \newline
{\it Case 2a:} \newline
$~~~~$ If $a_4\neq 0$ then by scaling $Y$, we can assume that
$a_4=1$. Now by the action of ${\rm Ad}(\exp a_2Y_2))$ on $Y$, we
can cancel the coefficient of $Y_2$:
\begin{eqnarray}
\overline{Y}=a_1Y_1+a_3Y_3+Y_4.
\end{eqnarray}
Then by applying ${\rm Ad}(\exp(a_3Y_3))$ on $\overline{Y}$ the
coefficient of $Y_3$ can be vanished and we have
\begin{eqnarray}
\overline{Y}'=a_1Y_1+Y_4.
\end{eqnarray}
The one-dimensional subalgebra generated by $Y$ is equivalent to
one generated by $\pm Y_1+Y_4$ which introduce parts {\it 8)} and
{\it 9)} of the theorem.\newline
{\it Case 2b:} \newline
$~~~~$ Let $a_4=0$ then  $Y$ is in the
form
\begin{eqnarray}
\widehat{Y}=a_1Y_1+a_2Y_2+a_3Y_3.
\end{eqnarray}
Suppose that $a_2\neq 0$ then if necessary we can let it equal to
$1$ and we obtain
\begin{eqnarray}
\widehat{Y}'=a_1Y_1+Y_2+a_3Y_3.
\end{eqnarray}
By acting ${\rm Ad}(\exp(a_3Y_3))$ on $\widehat{Y}'$, it changed
to $a_1Y_1+Y_2$ :\newline
{\it Case 2b-1:} \newline
$~~~~$ Let $a_1$ be nonzero. In this case we can make the
coefficient of $Y_1$ in $\widehat{Y}$ either $+1$ or $-1$ and find
{\it 6)}, {\it 7)} sections of the theorem.\newline
{\it Case 2b-2:} \newline
$~~~~$ If $a_1$ is zero then $Y_2$ is remained. Hence this case
suggests part {\it 2)}.\newline
{\it Case 2c:} \newline
$~~~~$ Finally if in the latter case $a_2$ be zero, then no
further simplification is possible and then $Y$ is one of
 cases of (\ref{eq:16}).
\medskip \noindent There is not any more possible case for studying and the
proof is complete.\hfill\ $\Box$

The coefficients $f,g$ of Eq. (\ref{eq:2}) depend on the variables
$x,u$. Therefore, we take their optimal system's projections on
the space $(x,u,f,g)$. The nonzero in $x$-axis or $u$-axis
projections of (\ref{eq:16}) are
\begin{eqnarray}\hspace{-0.7cm}\begin{array}{rlrl}
1)&Z^1=Y^2=Y^6=Y^7={\p}_x, \hspace{1mm}  &4)&Z^4=Y^5=Y^{10}=Y^{11}={\p}_x+2f{\p}_f+g{\p}_g,  \\
2)&Z^2=Y^3={\p}_u, \hspace{1mm}  &5)&Z^5=Y^{12}=Y^{14}=Y^{15}={\p}_x+u{\p}_u, \\
3)&Z^3=Y^4=Y^8=Y^9=-Y^{16}=u{\p}_u-2f{\p}_f-g{\p}_g, \hspace{1mm}
&6)&Z^6=Y^{13}=Y^{17}={\p}_x-u{\p}_u+4f{\p}_f+2g{\p}_g,
\label{eq:17}
\end{array}\end{eqnarray}
\begin{table}
\centering{\caption{The result of the
classification}}\label{table:3} \vspace{-0.35cm}\begin{eqnarray*}
\hspace{-0.75cm}\begin{array}{l l l l l l} \hline
  N       &\hspace{1cm} Z     &\hspace{1.1cm} \mbox{Invariant}  &\hspace{1cm} \mbox{Equation}
  &\hspace{1cm} \mbox{Additional operator}\,X^{(2)} \\ \hline
  1       &\hspace{1cm} Z^1   &\hspace{1.1cm} u         &\hspace{1cm} u_t=\Phi u_x^2+\Psi u_{xx}
  &\hspace{1cm} {\p}_x, {\p}_t+{\p}_x,\;-{\p}_t+{\p}_x \\
  2       &\hspace{1cm} Z^2   &\hspace{1.1cm} x         &\hspace{1cm} u_t=\Phi u_x^2+\Psi u_{xx}
          &\hspace{1cm} {\p}_u \\
 3        &\hspace{1cm} Z^3   &\hspace{1.1cm} x         &\hspace{1cm} u_t=u^2\Phi u_x^2+u\Psi u_{xx}
  &\hspace{1cm} t{\p}_t+u{\p}_u,\;(t+1){\p}_t+u{\p}_u,\;(t-1){\p}_t+u{\p}_u \\
 4        &\hspace{1cm} Z^4    &\hspace{1.1cm} u         &\hspace{1cm} u_t=e^{x^2}\Phi u_x^2+e^x\Psi u_{xx}
  &\hspace{1cm} {\p}_x, {\p}_t+{\p}_x,\;-{\p}_t+{\p}_x \\
 5        &\hspace{1cm} Z^5    &\hspace{1.1cm} {u\over{e^x}}&\hspace{1cm} u_t=\Phi u_x^2+u\Psi u_{xx}
          &\hspace{1cm} t{\p}_t+{\p}_x+u{\p}_u,\;(t+1){\p}_t+{\p}_x+u{\p}_u \\
 6        &\hspace{1cm} Z^6   &\hspace{1.1cm} -{1\over u}&\hspace{1cm} u_t=e^{x^4}\Phi u_x^2+e^{x^2}\Psi u_{xx}
  &\hspace{1cm} -t{\p}_t+{\p}_x-u{\p}_u,\;-(1+t){\p}_t+{\p}_x-u{\p}_u \\
  \hline
\end{array}\end{eqnarray*}
\end{table}
$~~~${\bf Proposition 4.2.} {\it Let ${\goth g}_m:=\langle Y_1,
\ldots, Y_m\rangle$, be an $m$-dimensional algebra. Denote by $Y^i
(i=1, \ldots, r, 0<r\leq m, r\in{\Bbb N})$ an optimal system of
one-dimensional subalgebras of ${\goth g}_m$ and by $Z^i\, (i =
1,\cdots, t, 0<t\leq r, t\in{\Bbb N})$ the projections of $Y^i$,
i.e., $Z^i = {\rm pr}(Y^i)$. If equations
\begin{eqnarray}
f = \Phi(x,u),\hspace{0.75cm} g =\Psi(x,u),\label{eq:18}
\end{eqnarray}
are invariant with respect to the optimal system $Z^i$ then the
equation
\begin{eqnarray}
u_t = \Phi(x,u)u_x^2+\Psi(x,u)u_{xx},\label{eq:19}
\end{eqnarray}
admits the operators $X^i=$ projection of $Y^i$ on $(t, x, u)$.}

$~~~${\bf Proposition 4.3.} {\it Let Eq. (\ref{eq:19}) and the
equation
\begin{eqnarray}
u_t = \Phi'(x,u)u_x^2+\Psi'(x,u)u_{xx},\label{eq:20}
\end{eqnarray}
be constructed according to Proposition 4.2. via optimal systems
$Z^i$ and ${Z^i}'$, respectively. If the subalgebras spanned on
the optimal systems $Z^i$ and ${Z^i}'$, respectively, are similar
in ${\goth g}_m$, then Eqs. (\ref{eq:19}) and (\ref{eq:20}) are
equivalent with respect to the equivalence group $G_m$, generated
by ${\goth g}_m$. }

$~~~$ Now we apply Proposition 4.2. and  Proposition 4.3. to the
optimal system (\ref{eq:17}) and obtain all nonequivalent Eq.
(\ref{eq:2}) admitting $\EE$-extensions of the principal Lie
algebra ${\goth g}$, by one dimension, i.e., equations of the form
(\ref{eq:2}) such that they admit, together with the one basic
operators (\ref{eq:21}) of ${\goth g}$, also a second operator
$X^{(2)}$. For every case, when this extension occurs, we indicate
the corresponding coefficients $f, g$ and the additional operator
$X^{(2)}$.

$~~~~$ We perform the algorithm passing from operators
$Z^i\,(i=1,\cdots,6)$ to $f, g$ and $X^{(2)}$ via the following
example. \newline
$~~~~$ Let consider the vector field
\begin{eqnarray}
Z^6={\p}_x-u{\p}_u+4f{\p}_f+2g{\p}_g,\label{eq:21}
\end{eqnarray}
then the characteristic equation corresponding to $Z^6$ is
\begin{eqnarray}
dx=\frac{du}{-u}=\frac{df}{4f}=\frac{dg}{2g},
\end{eqnarray}
and can be taken in the form
\begin{eqnarray}
I_1=ue^x,\hspace{5mm}I_2=\frac{f}{e^{x^4}},\hspace{5mm}I_3=\frac{g}{e^{x^2}}.
\end{eqnarray}
From the invariance equations we can write
\begin{eqnarray}
I_2=\Phi(I_1), \hspace{1cm} I_3= \Psi(I_1),
\end{eqnarray}
it follows that
\begin{eqnarray}
f=e^{x^4}\Phi(\lambda),\hspace{1cm} g=e^{x^2}\Phi(\lambda),
\end{eqnarray}
where $\lambda=I_1$.

$~~~$From Proposition 4.2. applied to the operator $Z^6$ we obtain the additional operator $X^{(2)}$
\begin{eqnarray}
-t{\p}_t+{\p}_x-u{\p}_u,\;\;\;-(1+t){\p}_t+{\p}_x-u{\p}_u.
\end{eqnarray}
After similar calculations applied to all operators (\ref{eq:17})
we obtain the following result (Table 3) for the preliminary group
classification of Eq. (\ref{eq:2}) admitting an extension ${\goth
g}_3$ of the principal Lie algebra ${\goth g}_1$.
\section{Conclusion}
In this paper, following the classical Lie method, the preliminary
group classification for the class of general Burgers' equation
(\ref{eq:2}) and investigated the algebraic structure of the
symmetry groups for this equation, is obtained. The classification
is obtained by constructing an optimal system with the aid of
Propositions 4.2. and 4.3.. The result of the work is summarized
in Table 3. Of course it is also possible to obtain the
corresponding reduced equations for all the cases in the
classification reported in Table 3. We omitted these for brevity.
%

\end{document}